\documentclass[a4paper,12pt]{amsart}\usepackage{amsfonts,amssymb,oldgerm,enumerate,mathrsfs,bbm,hyperref,rotating}
\UseRawInputEncoding

\def\C{\mathbb{C}}
\def\k{\mathbbm{k}}
\def\N{\mathbb{N}}
\def\bS{\mathbb{S}}\def\R{\mathbb{R}}\def\Z{\mathbb{Z}}

\def\bC{{\bar{C}}}

\def\di{\partial}
\def\bl{\langle}\def\br{\rangle}

\newcommand{\quot}[2]{{\left.\raisebox{1.6ex}{\footnotesize$#1$}  \!\!\!\!{\ensuremath\diagup}\!\!\raisebox{-1ex}{\footnotesize$#2$}\right.}}
\renewcommand{\stackrel}[2]{\ \lower 0.2ex \hbox{$\mathrel{\mathop{#2}\limits^{#1}}$}\ }

\makeatletter
\newcommand{\tpitchfork}{%
  \vbox{
    \baselineskip\z@skip
    \lineskip-.52ex
    \lineskiplimit\maxdimen
    \m@th
    \ialign{##\crcr\hidewidth\smash{$-$}\hidewidth\crcr$\pitchfork$\crcr}
  }%
}
\makeatother

\def\tga{{\tilde{\ga}}} \def\tbk{{\tilde{\k}}}
\def\tt{\tilde{t}}
\def\cU{\mathcal U}\def\tu{\tilde{u}}
\def\tU{\tilde{U}}\def\tcU{\tilde{\cU}}

\def\hx{\hat{x}}\def\hy{\hat{y}}

\def\ga{\gamma}\def\de{\delta}
\def\ep{\epsilon}\def\om{\omega}
\def\si{\sigma}\def\Si{\Sigma}\def\tSi{\tilde\Si}

\def\cC{\mathscr C}
\def\cG{\mathcal G}\def\cK{{\mathcal K}}\def\cR{\mathscr{R}}

\def\uf{{\underline{f}}}\def\uh{{\underline{h}}}

 \def\uom{{\underline{\om}}}\def\up{{\underline{p}}}\def\us{{\underline{s}}}

\def\ux{{\underline{x}}}\def\uz{{\underline{z}}}

\def\one{{1\hspace{-0.1cm}\rm I}}
\def\empty{\varnothing}

\newcommand{\bbm}{\begin{bmatrix}}\newcommand{\ebm}{\end{bmatrix}}
\newcommand{\ber}{\begin{array}{l}}\newcommand{\eer}{\end{array}}
\newcommand{\bpm}{\begin{pmatrix}}\newcommand{\epm}{\end{pmatrix}}
\newcommand{\bM}{\begin{matrix}}\newcommand{\eM}{\end{matrix}}
\newcommand{\bee}{\begin{enumerate}}\newcommand{\eee}{\end{enumerate}}
\newcommand{\bei}{\begin{itemize}}\newcommand{\eei}{\end{itemize}}

\def\iff{if and only if }
\def\sset{\subset}\def\sseteq{\subseteq}\def\smin{\setminus}

\newcommand{\beq}{\begin{equation}}\newcommand{\eeq}{\end{equation}}  %\vspace{-0.1cm}{\vspace{-0.2cm}
\newtheorem{Lemma}{Lemma}[section]\newcommand{\bel}{\begin{Lemma}}\newcommand{\eel}{\end{Lemma}}
\newtheorem{Example}[Lemma]{Example}\newcommand{\bex}{\begin{Example}\rm}%\newcommand{\eex}{\hfill$\lrcorner$\end{Example}}
\newcommand{\eex}{\end{Example}}
\newtheorem{Proposition}[Lemma]{Proposition}\newcommand{\bprop}{\begin{Proposition}}\newcommand{\eprop}{\end{Proposition}}
\newtheorem{Property}[Lemma]{Property}\newcommand{\bproperty}{\begin{Property}}\newcommand{\eproperty}{\end{Property}}
\newtheorem{Definition-Proposition}[Lemma]{Definition-Proposition}

\def\bpr{~\\{\em Proof.\ }}
\newcommand{\epr}{{\hfill\ensuremath\blacksquare}\\}
\newtheorem{Theorem}[Lemma]{Theorem}\newcommand{\bthe}{\begin{Theorem}}\newcommand{\ethe}{\end{Theorem}}
\newtheorem{Definition}[Lemma]{Definition}\newcommand{\bed}{\begin{Definition}}\newcommand{\eed}{\end{Definition}}
\newtheorem{Remark}[Lemma]{Remark}\newcommand{\beR}{\begin{Remark}\rm}\newcommand{\eeR}{\end{Remark}}
\newtheorem{Corollary}[Lemma]{Corollary}\newcommand{\bcor}{\begin{Corollary}}\newcommand{\ecor}{\end{Corollary}}
\newcommand{\bet}{\begin{tabular}{cccccccc}}\newcommand{\eet}{\end{tabular}}

\title[]{D\MakeLowercase{eforming the weighted-homogeneous foliation, and trivializing families of semi-weighted homogeneous} ICIS.}
\author[]{D\MakeLowercase{mitry}  K\MakeLowercase{erner and} R\MakeLowercase{odrigo} M\MakeLowercase{endes}}

\thanks{We were supported by the Israel Science Foundation,  grants No.  1910/18 and 1405/22}

\address{Dmitry Kerner: Department of Mathematics, Ben Gurion University of the Negev, P.O.B. 653, Be'er Sheva 84105, Israel. dmitry.kerner@gmail.com}
\address{Rodrigo Mendes: Instituto de ci\^encias exatas e da natureza, Universidade de Integra\c{c}\~ao Internacional da Lusofonia Afro-Brasileira (unilab), Campus dos Palmares, Cep. 62785-000. Acarape-Ce, Brasil and  Departament of Mathematics, Ben Gurion University of the Negev, P.O.B. 653, Be'er Sheva 84105, Israel. rodrigomendes@unilab.edu.br}


\date{\today\ \  filename: \jobname.tex}

\begin{document}\setcounter{secnumdepth}{6} \setcounter{tocdepth}{1}

 \begin{abstract}
Let $X_o:=V(\uf_\up)\sset(\k^N,o)$ be a weighted-homogeneous complete intersection germ (with arbitrary singularities, possibly non-reduced).
 Let $\{\ga_\us\}_{\us\in \bS}$ be the corresponding foliation of $\k^N$ by weighted-homogeneous real arcs.
 Take a deformation by higher order terms, $X_\ep:=V(\uf_\up+\ep\cdot \uf_{>\up}).$ Does the foliation $\{\ga_\us\}_\us$ deform compatibly with $X_\ep?$ We identify the ``obstruction locus", $\Si\sset X_o,$ outside of which such a deformation does exist, and possesses exceptionally nice properties.

Using this deformed foliation we construct a contact  trivialization of the family $\uf_\up+\ep\cdot \uf_{>\up} $ by a homeomorphism  that is
  real analytic (resp. Nash) off the origin, differentiable at the origin, whose presentation in weighted-polar coordinates is globally real-analytic (resp. globally Nash), and with controlled Lipschitz/$C^1$-properties.
 \end{abstract}
 \maketitle
%\tableofcontents

\vspace{-0.1cm}
\section{Introduction}
\vspace{-0.1cm}

Let $\k=\R,\C.$ Fix some  integral positive weights of coordinates in $\k^N,$ i.e. $\om_1\le \om_2\le\cdots\le\om_N,$ $\om_i\in \Z_{>0}.$
Take a regular sequence $\uf_\up=(f_{p_1},\dots,f_{p_c})$ of   weighted-homogeneous polynomials, of
 weights $\up=(p_1,\dots,p_c),$ where   $p_i=ord_\om (f_{p_i}).$
 These polynomials define the
  weighted-homogeneous (complete intersection) germ $X_o:=V(\uf_\up)\sset (\k^N,o).$
   The germ $X_o$ can have arbitrary singularity, can be non-reduced.

 Take the weighted homogeneous arcs, $\ga_\us(t)=(s_1 t^{\om_1},\dots, s_N t^{\om_N}).$  Here $t\in \R_{\ge0}$ and $\us\in \bS\sset \k^N,$ the unit sphere.  These arcs foliate the space $\k^N,$ and this foliation is compatible with the germ $X_o\sset (\k^N,o).$
 This foliation is the first basic tool in the study of Geometry and Topology of such germs.

\subsection{}\label{Sec.Intro.1}   Perturb $X_o$ by higher-order terms, $X_\ep:=V(\uf_\up+\ep\cdot\uf_{>\up})\sset (\k^N,o),$ here $f_{>p_i}\in \k\{x\}$ and $ord_\om f_{p_i} < ord_\om f_{>p_i} .$
  One would like to deform the foliation $\{\ga_\us\}_\us$ to a foliation $\{\ga_{\ep,\us}\}_\us$ of $(\k^N,o),$  compatibly with $X_\ep.$
 This is not always possible, see Example \ref{Ex.Briancon.Speder-type}.  Such a deformation is in general obstructed by the locus $Sing(X_o)$ and by the singularities of the intersection of $X_o$ with the flag of coordinate planes, $\k^N\supset V(x_1)\supset\cdots\supset V(x_1,\dots,x_N).$ Altogether they give  ``obstruction locus" $\Si\sset X_o\sset (\k^N,o)$, defined in \S\ref{Sec.Obstruction.Locus}.
%\beq
%  \Si =Sing[X_o] \cup  Sing[X_o\cap V(x_1)]\cup \cdots\cup  Sing[X_o\cap V(x_1,\dots,x_{N-1})]\sset (\k^N,o).
%  \eeq
\\{\bf Lemma \ref{Thm.semi-weighted.hom.foliation}} tells (roughly): {\em outside of a `hornic' neighborhood of this obstruction locus, $\cU(\Si)\sset (\k^N,o),$
 the foliation $\{\ga_\us\}_\us$ deforms to a foliation $\{\ga_{\ep,\us}\}_\us,$ compatibly with the family $X_\ep,$ preserving the (non)tangency of the arcs, and with various other controlled properties.}

In particular, if $\Si=\{o\},$ then $\{\ga_\us\}_\us$ deforms to a foliation $\{\ga_{\ep,\us}\}_\us$ on the whole $(\k^N,o).$

 This construction is purely algebraic, no vector field integration is involved. It is based on the implicit function theorem of Bourbaki-Tougeron.
  The functions $\ga_{\ep,\us}(t)$ are $\R$-analytic (in $\ep,\us,t$). If the perturbations $\uf_{>\up}$ are Nash (e.g. polynomials), then  $\ga_{\ep,\us}(t)$  are Nash functions.
 %(See also [Kerner-Mendes.0]).
\\The deformed ``semi-weighted homogeneous" foliation  $\{\!\ga_{\ep,\us}\!\}_\us$ is as  useful for $X_\ep,$  as the weighted homogeneous foliation $\{\!\ga_\us\!\}_\us$  is for $X_o.$

\subsection{} A related question is the trivialization of the family $X_\ep\!:=\!V(\uf_\up\!\!+\!\ep\!\cdot\!\uf_{>\up}).$ Suppose the (real or complex) singular set of $X_o$ is just one point. Then the family of germs $X_\ep$ is embedded topologically trivial, while the family of functions
 $\uf_\up\!\!+\!\ep\!\cdot\!\uf_{>\up}$ is contact  topologically trivial.
 Various stronger trivializations were obtained in the last 70 years (see example \ref{Ex.Trivialiations} for  details).

 \bee[\!\!$\bullet$]
 \item
 The contact-Lipschitz/$C^k$ trivializations were given (under additional assumptions) in  \cite{Ruas-Saia.97},  \cite{Fernandes-Ruas},  \cite{Costa-Saia-Soares}.
 The construction went by integrating Lipschitz/$C^k$ vector fields. The resulting trivialization was not subanalytic/semialgebraic.

 \item
 The family $f_p+\ep\cdot f_{>p}$ (where the weighted-homogeneous polynomial $f_p$ has an isolated critical point) was right-trivialized in \cite{Fukui.Paunescu.00}
 by a homeomorphism whose strict transform  under a weighted-homogeneous blowup is real-analytic.
 (The homogeneous case was done in \cite{Kuo.80}.)
 That trivialization does not seem to have controllable Lipschitz/$C^1$-properties.
 \eee

\medskip

\noindent The deformed foliation $ \ga_{\ep,\us} $ of \S\ref{Sec.Intro.1} helps to trivialize the family $\uf_\up\!+\!\ep\!\cdot\!\uf_{>\up} $ in a very special way.
 \\{\bf  Theorem \ref{Thm.Trivialization.of.Def.by.higher.order.terms}} (roughly): {\em If $\Si=\{o\}$ then the family $\uf_\up\!+\!\ep\!\cdot\!\uf_{>\up} $ is contact-trivializable by a homeomorphism that is $\R$-analytic on $(\k^N,o)\smin o,$   differentiable at $o,$ becomes globally $\R$-analytic in the  weighted-polar coordinates, and has controlled Lipschitz/$C^1$-properties.
\\Moreover, if the perturbations $\uf_{>\up}$ are Nash functions (e.g. polynomials), then the trivialization is Nash  on $(\k^N,o)\smin o$ and  becomes globally Nash in the weighted-polar coordinates.}

\beR
\bee[\!\!\!\bf i.\!]
\item All the previous trivializations were obtained by vector field integration. This left no hope for Nash trivialization.

Recall that   analyticity in weighted-polar coordinates is a much stronger property than subanalyticity.
\item
We can allow perturbations of the same order,   $ord_\om f_{p_i} \le ord_\om f_{\ge p_i},$ but then one must assume $|\ep|\ll1,$ and not all the properties of Lemma \ref{Thm.semi-weighted.hom.foliation}, Theorem \ref{Thm.Trivialization.of.Def.by.higher.order.terms} hold. See Remark \ref{Rem.Deformation.of.Weight.Hom.Foliation.by.other.terms.} and Corollary \ref{Thm.Trivialization.for.defs.by.same.order.terms}.
\item While in this paper we work only with $\k=\R,\C,$ such trivialization/deformation questions are important over arbitrary fields, of any characteristic, \cite{B.G.K.}. An advantage of our method is: if $\uf_\up,\uf_{>\up}$ are defined over a field $\tbk,$ then the deformed foliation (of Lemma \ref{Thm.semi-weighted.hom.foliation}) is defined over the same $\tbk.$ (And similarly for the trivialization of Theorem \ref{Thm.Trivialization.of.Def.by.higher.order.terms}.)
\eee
\eeR
\noindent Having such a trivialization of the family $X_\ep$ (resp. the deformed foliation $\ga_{\ep,\us}$) is extremely useful for the study of semi-weighted homogeneous germs, Newton-non-degenerate germs, and other ``perturbations of weighted-homogeneous germs".
 In \cite{Kerner-Mendes.Weighted.Homogen} we apply this deformed foliation to detect fast vanishing cycles on perturbations of weighted-homogeneous germs. More generally, this trivialization allows to extend several recent results of Lipschitz Geometry of Singularities from weighted-homogeneous germs to broader classes of germs.

\subsection{Acknowledgements} Our thanks are to Lev Birbrair, Alexandre Fernandes, Andrey Gabrielov, Edson Sampaio, for important advices.

\subsection{Notations and conventions}\label{Sec.Notations}  By   germs we mean their small representatives.
\bee[\!\!\bf i.\!]

\item
Take the standard sphere   $\bS:=\{ x|\ \| x\|=1\}\sset \k^N.$  It is $ S^{N-1}\sset \R^N$ or $ S^{2N-1}\sset \C^N.$
\\
We use multi-indices, $\ux=(x_1,\dots,x_N)\in \k^N,$  $\us=(s_1,\dots,s_N)\sset \bS,$ $\uom=(\om_1,\dots,\om_N)\in \Z^N_{>0},$ $\up=(p_1,\dots,p_c),$ $\uf_\up=(f_{p_1},\dots,f_{p_c}).$
 Accordingly  we abbreviate:
 \beq
 t^{\uom}\cdot \us\!: =\!(t^{\om_1}  s_1 ,\dots,t^{\om_N} s_N  ), \ \  \
 t^{\up}\cdot \uf_\up\!: =\!(t^{p_1}  f_{p_1} ,\dots,t^{p_c}  f_{p_c}  ), \  \ \
  t^{-\up}\cdot \uf_{>\up}\!: =\!(t^{-p_1}  f_{\!>p_1} ,\dots,t^{-p_c}  f_{\!>p_c}  ).
 \eeq%\vspace{-0.1cm}
We always assume: $gcd(\om_1,\dots,\om_N)=1.$

Occasionally we pass to the  weighted polar coordinates, $\bS\times\R_{\ge0}\stackrel{\si}{\to}\k^N,$   $(\us,t)\to t^\uom\cdot \us.$
 For an analytic subgerm $\Si\sset (\k^N,o)$ we get ``the strict transform", $\tSi=\overline{\si^{-1}(\Si\smin o)}\sset \bS\times\R_{\ge0}.$
 We get also ``the link at $o$", $Link_o[\Si]:=\tSi\cap \{t=o\}.$

\item All our arcs are real-analytic.

A germ $X$ is called ``foliated by a set of arcs $\{\ga_\us\}_\us$", if there exists a representative of $X$ that is covered by those arcs, and moreover (for this representative): $\ga_\us\cap \ga_{\us'}=o$ for $\us\neq\us'.$

  The tangency order of  arcs $\ga_1, \ga_2 \sset (k^N,o)$ is $tord_{(k^N,o)} (\ga_1,\ga_2):=ord_t \|\ga_1(t)-\ga_2(t)\|.$
 Here we take the  length-parametrization, thus $\ga(t)\!=\!(s_1\cdot t\!+\!o(t),\dots,s_N\cdot  t\!+\!o(t)),$ with $\us\!\in\!\bS.$

\item
 Take a   weighted-homogeneous complete intersection $X_o:=V(\uf_\up)\sset (\k^N,o),$ here
  $\uf_\up=(f_{p_1},\dots,f_{p_c})$ is a regular sequence of polynomials,
  and $ord_\om f_{p_i}=p_i,$ i.e. $f_{p_i}(t^{\om_1} x_1,\dots,t^{\om_N}x_N)=t^{p_i}\cdot f_{p_i}(x_1,\dots,x_N).$ W.l.o.g. we assume $(\uf_\up)\sseteq (\ux)^2\sset \k[\ux],$ thus $X_o$ is singular.  The germ $X_o$  can be non-reduced.
\\
The singular locus $Sing(X_o)$ is defined (inside $X_o$) by the determinantal ideal of the Jacobian matrix $I_c[f'_1,\dots,f'_c].$ Here $\{f'_j\}$ are the column-gradients.  Thus $Sing(X_o)$ is a weighted-homogeneous $\k$-algebraic subgerm.
The dimension $dim_\k Sing(X_o)$ is taken set-theoretically. For $\k\!=\!\C$ this dimension can be computed also via the height of the defining ideal.

\eee

\vspace{-0.2cm}
\section{Deforming the weighted-homogeneous foliation $\{\ga_\us\}_\us$}\label{Sec.Foliations.of.semi-Weight.Hom.germs}
\vspace{-0.1cm}
The arcs
 $\ga_\us(t)=t^{\underline{\om}}\cdot\us:=(t^{\om_1} s_1,\dots,t^{\om_N}s_N)$ with $\us\in \bS, $
   form the weighted-homogeneous foliation of $\k^N.$ It is compatible with $X_o\sset \k^N,$ i.e. either $\ga_s\sset X_o$ or $\ga_s\cap X_o=\{o\}.$

 Take a deformation by higher order terms $X_\ep:=V(\uf_{\up}+\ep\cdot \uf_{>\up})\sset (\k^N,o),$  for $\ep\in Disc\sset\k^1 ,$ i.e. $|\ep|\le1.$
  Here $ord_\om(f_{p_i})<ord_\om (f_{>p_i}).$   This deformation is flat, as $X_o$ is a complete intersection.
  This family is not equisingular in whichever sense, as we did not assume that $X_o$ has an isolated singularity.

One would like to deform the foliation $\{\ga_\us\}$ to a foliation $\{\ga_{\ep,\us}\}$ of $(\k^N,o)$ (by terms of higher weights),
    compatibly with $X_\ep,$ i.e.: $\ga_\us\sset X_o$ \iff $\ga_{\ep,\us}\sset X_\ep.$
          The minimal condition to impose is ``each arc is tangent to its deformation", i.e.:
    \beq
    tord_{(\k^N,o)}(\ga_{ \us},\ga_{\ep,\us})>1 \text{ for each }\us\in Link[X_o]=\bS \cap X_o \text{ and }\ep\in Disc.
     \eeq
For $\k\!=\!\C$    such a deformation exists, e.g., if $X_o$ has an isolated singularity, and the family $\{X_\ep\}_\ep$ is equi-resolvable by a weighted blowup, and with additional assumptions on the exceptional divisor.
Just the $\mu\!=\!const$ assumption  on family $\{X_\ep\}_\ep$  is  not sufficient,  due to   examples of  Brian\c{c}on-Speder-type.

\bex\label{Ex.Briancon.Speder-type}
Let $f_\ep(x,y,z)=x^{p_x}+z^{p_z}+xy^{p_y}+\ep\cdot y^{p_y+1},$ with the weights $(\om_x,\om_y,\om_z)= (\frac{1}{p_x},\frac{1-\frac{1}{p_x}}{p_y},\frac{1}{p_z}).$ (Up to $\Z$-scaling.)
 Thus $V(f_o)\sset (\C^3,o)$ has a weighted-homogeneous isolated singularity.
 The perturbation $\ep\cdot y^{p_y+1}$ is of weight$> ord_\om f_o$ \iff $\frac{p_y+1}{p_y}(1-\frac{1}{p_x})>1,$ i.e. $p_x>p_y+1.$
 In this case $f_\ep$ is a $\mu=const$ family. Take the arc $\ga_o(t)=(0,t,0)\sset V(f_o).$
 Its deformation by higher  order terms (to ensure $tord(\ga_o,\ga_\ep)>1$) can be presented as $\ga_\ep(t)=(t^{d_x}(\dots),t,t^{d_z}(\dots)),$ for some $d_x,d_z>1.$ The condition $\ga_\ep\sset V(f_\ep)$ is then:
  $t^{p_x d_x}(\dots)+t^{p_z d_z}(\dots)+t^{d_x+p_y}(\dots)+\ep t^{p_y+1}=0.$
Assuming $p_z>p_y+1$ we have: $p_x d_x,p_z d_z,d_x+p_y>p_y+1.$ And this contradicts the condition $f_\ep(\ga_\ep(t))=0.$
\eex

In this example  ``the problematic arc" $\ga_o$ lies inside the intersection $X_o\cap V(x,z).$  In fact,  $X_o\cap V(x,z)\sset \k^3$ is
 a line, while the expected dimension of this intersection is zero.
 We show that such intersections of unexpected dimension, or with unexpected singularities, are {\em the only} obstructions to   deforming the foliation $\{\ga_\us\}_\us$ to $\{\ga_{\ep,\us}\}_\us.$

\subsection{The obstruction locus $\Si$}\label{Sec.Obstruction.Locus}
Some of the weights $\om_\bullet$ might coincide, thus we split them:\vspace{-0.1cm}
\beq\label{Eq.Splitting.of.weights}
\om_1=\cdots=\om_{r_1}<\om_{r_1+1}=\cdots=\om_{r_2}<\cdots<\om_{r_k+1}=\cdots=\om_N.
\eeq
  Intersect $X_o$ by the coordinate planes $V(x_1,\dots,x_{r_j})\sset (\k^N,o)$ and take the union of the corresponding   singular loci (see \S\ref{Sec.Notations}):
  \beq\label{Eq.Def.of.Sigma}
  \Si\!:=\!Sing[X_o]\! \cup \! Sing[X_o\!\cap\! V(x_1,\dots,x_{r_1})] \!\cup \! Sing[X_o\!\cap\! V(x_1,\dots,x_{r_2})]\! \cup\! \cdots
  \! \cup\!   Sing[X_o\!\cap\! V(x_1,\dots,x_{r_k})].
\eeq
  %\hspace{1cm}
% \vspace{-0.2cm} \vspace{-0.2cm}  \[\hspace{11cm}\vspace{-0.2cm} \cdots    \! \cup\!   Sing[X_o\!\cap\! V(x_1,\dots,x_{r_k})].  \]
Here each intersection   $X_o\!\cap\! V(x_1,\dots,x_{r_j})$ is considered as  a scheme with ``prescribed dimension":
\bei
\item (for $\k=\R$)  $dim_\R\le dim_\R(X_o)-r_j.$
\item (for $\k=\C$)  $dim_\C= dim_\C(X_o)-r_j.$
\eei
 Therefore $Sing[X_o\cap V(x_1,\dots,x_{r_j})]$ contains all the non-reduced components of $X_0\cap V(x_1,\dots,x_{r_j})$, and also all the components of
  dimension$>dim_\k(X_o)-r_j.$

Altogether,  $\Si  \sset (\k^N,o)$ is a weighted-homogeneous ($\k$-algebraic) subgerm.
\bex\label{Ex.Sigma}
\bee[\bf i.]
\item If $\om_1=\cdots=\om_N$ then $r_1=N$ and $\Si=Sing(X_o).$
\item Suppose $\om_1\!<\!\om_2\!<\!\cdots\!<\!\om_{n+1}$ and let $X_o\!\sset \!(\k^{n+1},o)$ be a weighted-homogeneous hypersurface germ with isolated singularity,      i.e. $Sing(X_o)\!=\!\{o\}.$ Suppose $dim_\k( X_o\cap V(x_1,\dots, x_j))\!\le\! n-j,$ as one could expect.
     Moreover, suppose all the germs  $X_o\!\cap \!V(x_1 ),$ \dots, $X_o\!\cap \!V(x_1,\dots, x_{n-1})$ are reduced outside of $o,$ and $Sing[X_o\!\cap\! V(x_1,\dots, x_j)]\!=\!\{o\}.$
      Then $\Si=\{o\}.$
\item (Continuing example \ref{Ex.Briancon.Speder-type}) Let $f_o(x,y,z)=x^{p_x}+z^{p_z}+xy^{p_y},$ suppose the weights satisfy: $\om_x<\om_z<\om_y,$ i.e. $p_x>p_z>p_y.$ Then $\Si=V(x,z)\sset (\k^3,o).$
\eee
\eex

\subsection{Deforming the foliation}
Having the obstruction locus, $\Si\sset (\k^N,o),$ the naive guess is: the arcs $\{\ga_\us\}_\us$ deform to the arcs $\{\ga_{\ep,\us}\}_\us,$
 and for each $\ep$ these later arcs foliate the germ $(\k^N,o)\smin \Si.$ This does not always hold. In fact, take a small representative $Ball\smin \Si$ of  $(\k^N,o)\smin \Si.$ It is not necessary covered by the arcs $\{\ga_{\ep,\us}\}_\us,$ as the radius of convergence of
  $\ga_{\ep,\us}(t)\in \R\{t\}$ can vanish when $\us$ approaches $Link[\Si].$
The correct statement will be: the arcs $\{\ga_{\ep,\us}\}_\us$ foliate the germ  $(\k^N,o)\smin \cU(\Si)$ for certain `hornic' neighborhood
 $\Si\sset \cU(\Si).$
\\Take the change of $\om$-orders,   $\de\!:=\!min_i\{ ord_\om f_{>p_i}\!-\! ord_\om f_{p_i} \}.$
  If $f_{>p_i}=0,$ one puts $ord_\om f_{>p_i}=\infty$.
 Thus $0<\de\in \N\cup\{\infty\}.$
 \bel\label{Thm.semi-weighted.hom.foliation}
There exists  a family of arcs $\{\ga_{\ep,\us}\}_\us$, for $\ep\in  Disc$ and  $\us\in \bS\smin Link[\Si]$,
\beq
\ga_{\ep,\us}(t):=t^{\uom}\cdot(\us+\ep\cdot t^\de\cdot \uh):=(t^{\om_1}  s_1+t^{\om_1+\de} \ep\cdot h_1(\ep,\us,t) ,\ \dots\ ,t^{\om_N} s_N+t^{\om_N+\de} \ep\cdot h_N(\ep,\us,t) ),
\quad
\eeq
with the following properties.

\bee[\bf \!\!\!1.\!]
\item
For each $\ep,\us$ the arcs  are  real-analytic, i.e.   $\ga_{\ep,\us}(t)\in \R\{t\}.$
  If the perturbations $\uf_{>\up}$ are Nash function-germs (e.g. polynomials), then the arcs   are  Nash, i.e.
   $\ga_{\ep,\us}(t)\in \R\bl t\br.$

 \item For all $\ep,\ep'\in Disc$ and $\us,\us'\in \bS \smin Link[\Si]:$
\bee[\!\!\bf i.]
\item  $\ga_\us\sset X_o$ \iff   $\ga_{\ep,\us}\sset X_\ep.$

 \item $tord_{\k^N}(\ga_\us,\ga_{\ep,\us})\ge 1+\frac{\de}{\om_N}.$
 \item   In particular, $tord_{\k^N}(\ga_\us,\ga_{ \us'})=1$ \iff $tord_{\k^N}(\ga_{\ep,\us},\ga_{\ep',\us'})=1.$
% \item If  $tord_{\k^N}(\ga_\us,\ga_{ \us'})>1$ then $tord_{\k^N}(\ga_{\ep,\us},\ga_{\ep',\us'})>1.$

\item    $\ga_\us\sset V(x_1,\dots,x_{r_j})$ \iff $\ga_{\ep,\us}\sset V(x_1,\dots,x_{r_j}).$
\item If $\ga_\us\in X_o\cap X_\ep\smin \Si$ (for some $\ep>0$ and $\us\in\bS$), then $\ga_{\ep,\us}=\ga_\us$ for each $\ep.$
 \eee
 \item
(The case $\Si=\{o\},$ thus $Link[\Si]=\empty$)
\bee[\!\!\bf i.]
\item The arcs are  analytic in  $t,\ep,\us,$ i.e. $\ga_{\ep,\us}(t)\in C^\om(Disc\times\bS\times[0,t_o)),$ for some $0<t_o.$
\\If the perturbations $\uf_{\!>\up}$ \!are Nash function-germs, \!then  $\ga_{\ep,\us}(t)\!\in\! Nash(Disc\!\times\!\bS\!\times\![0,t_o)).$
\item
For each  $\ep\in Disc$ the arcs $\ga_{\ep,\us}$ foliate the whole germ $(\k^N,o).$
\eee
\item (The case $\Si\neq\{o\},$ thus $Link[\Si]\neq\empty$) For any ($\R$-analytic) neighborhood $Link[\Si]\sset \tU\sset\bS$ there exists $0<t_o$
 and $o\in Ball_{r_o}\sset \k^N,$ satisfying:
\bee[\!\!\bf i.\!]
\item
$\ga_{\ep,\us}(t)\in C^\om(Disc\times(\bS\smin \tU)\times[0,t_o)).$
\\If the perturbations $\uf_{>\up}$ are Nash function-germs, then  $\ga_{\ep,\us}(t)\!\in\! Nash(Disc\!\times\!(\bS\smin \tU)\!\times\![0,t_o)).$

\item $\ga_{\ep,\us}\cap \ga_{\ep,\us'}\cap Ball_{r_o}=\{o\}$ for all $\us,\us'\in \bS\smin\tcU $ with $s\neq s'.$
\item There exists an $\R$-analytic neighborhood $\Si\sset \cU(\Si)\sset (\k^N,o)$ with $Link[  \cU(\Si) ]=\tU,$
  $Link_o[  \cU(\Si) ]=Link_o[\Si],$
 and such that (for each $\ep\in Disc$) the arcs cover the complement of $\cU(\Si),$
  i.e.
  \[(\cup_{\us}\ga_{\ep,\us})\supseteq ( \k^N \smin\cU(\Si) )\cap Ball_{r_o}.
  \]
\eee
\eee
\eel
\bpr
 We begin  with the foliation $\ga_{\ep,\us}(t)=t^{\underline{\om}}(\us+\ep\cdot t^\de\cdot \uh),$ where $\uh$ is the vector of unknowns.
\bee[\!\!\bf Step 1.]
\item
 The main condition on $\uh$  is: if  $\ga_\us\sset  X_o,$  i.e. $\uf_\up(\ga_\us(t))=0,$ then   $\ga_{\ep,\us}\sset X_\ep, $ i.e. $(\uf_\up+\ep  \cdot\uf_{>\up})(\ga_{\ep,\us}(t))=0.$ We will resolve this condition.
 Taylor-expand the first term:
\vspace{-0.5cm}
\beq
 \uf_\up(\ga_{\ep,\us}(t))=t^\up\cdot \uf_\up(\us+\ep\cdot t^\de\cdot \uh)=t^\up\cdot [\uf_\up(\us)\put(-25,-5){\vector(3,2){25}}\put(2,10){\scriptsize{0}}+\ep\cdot t^\de\cdot  \uf_\up'|_\us(\uh)+\frac{\ep^2  t^{2\de}}{2}\uf_\up''|_\us(\uh,\uh)+\dots].
 \eeq
Then the condition $\ga_{\ep,\us} \sset X_\ep$ becomes:
\beq\label{Eq.pre-IFT.form}
[\uf_\up'|_\us(\uh)+\frac{\ep  t^\de}{2}\cdot \uf_\up''|_\us(\uh,\uh)+\dots]+  t^{-\de}\cdot t^{-\up}\cdot \uf_{>\up}(t^{\underline{\om}}(\us+\ep\cdot \uh))=0,
\quad \quad \forall\ \us\in Link[X_o].
\eeq
In the chosen coordinates ($x_1\dots x_N$) the first derivative is the column of (rows of) gradients, $\uf_\up'|_\us=\nabla\uf_\up |_\us\in Mat_{c\times N}(\k[\us]).$
 We are looking for the solution of particular form, using the transposed gradients with rescaled entries\footnote{The factors $\sum |s^2_i|$ are inserted here to ensure   properties 2.iii,2.iv  and also Theorem \ref{Thm.Trivialization.of.Def.by.higher.order.terms}, as we will see later.},
\beq\label{Eq.inside.proof.Def.of.h}
\uh=\bigg[(\sum ^{r_1}_{i=1}|s^2_i|)\cdot[\di_{x_1}\uf_{\up}|_\us]^T\cdot \uz,\dots, (\sum^{r_1}_{i=1}|s^2_i|)\cdot[\di_{x_{r_1}}\uf_{\up}|_\us]^T\cdot \uz,
(\sum^{r_2}_{i=1}|s^2_i|)\cdot[\di_{x_{r_1+1}}\uf_{\up}|_\us]^T\cdot \uz,
\dots,\vspace{-0.2cm}
\eeq
\[
\hspace{7cm}
\dots,(\sum^{r_2}_{i=1}|s^2_i|)\cdot[\di_{x_{r_2}}\uf_{\up}|_\us]^T\cdot \uz,\dots
 (\sum^N_{i=1} |s^2_i|)\cdot[\di_{x_N}\uf_{\up}|_\us]^T\cdot \uz\bigg].
\]

 Here $\uz=( z_1,\dots, z_c)$ is a  $c$-column of unknowns.
 In the case $\k=\C$ one takes the complex-conjugates, $\overline{[\di_{x_j}\uf_{\up}|_\us]^T}\in Mat_{1\times r}(\C).$

  Condition \eqref{Eq.pre-IFT.form} turns (with this Ansatz) into:
 \beq\label{Eq.in.the.proof}
  \nabla\uf_\up|_\us \cdot \bbm (\sum ^{r_1}_{i=1}|s^2_i|)\cdot[\di_{x_1}\uf_\up|_\us]^T\\\dots\\(\sum^N_{i=1} |s^2_i|)\cdot[\di_{x_N}\uf_\up|_\us]^T\ebm  \cdot \uz
   +   t^{-\de}\cdot t^{-\up}\cdot\uf_{>\up}(t^\uom\cdot \us)   +   \ep\cdot H(\uz,\us,\ep,t) =0.
 \eeq
The ingredients are $\R$-analytic:   $ t^{-\de}\cdot  t^{-\up}\cdot\uf_{>\up}(t^\uom\cdot \us) \in C^\om(\bS\times[0,t_o)),$ and
   $H(\uz,\us,\ep,t)\in (\uz)^2+(t^{ \de}\cdot\uz)\in C^\om(Ball_\uz\times \bS\times Disc\times [0,t_o)).$

The first term can be presented (using the diagonal $N\times N$ matrix) as
\beq
\nabla \uf_\up|_\us \cdot diag\Big[\sum^{r_1}|s^2_i|,\dots,\sum^{r_1}|s^2_i|,\sum^{r_2}|s^2_i|,\dots,\sum^{r_2}|s^2_i| ,\dots, \sum^N |s^2_i|\Big]_{N\times N}\cdot [\nabla \uf_\up|_\us]^T.
\eeq
Accordingly, we introduce the matrix with rescaled gradients  (and semi-algebraic entries),
\beq
\tilde\nabla \uf_\up|_\us:=\Big[\sqrt{\sum ^{r_1}_{i=1}|s^2_i|}\cdot \di_{x_1} \uf_\up|_\us,\dots,\sqrt{\sum ^{r_1}_{i=1}|s^2_i|}\cdot\di_{x_{r_1}} \uf_\up|_\us,\dots,\sqrt{\sum^N_{i=1} |s^2_i|}\cdot\di_{x_N} \uf_\up|_\us\Big]\in Mat_{r\times N}.
\eeq
   Then  equation \eqref{Eq.in.the.proof} becomes:  $\tilde\nabla \uf_\up|_\us\cdot[\tilde\nabla \uf_\up|_\us]^T \cdot \uz
   + \ep\cdot t^{-\de}\cdot  t^{-\up}\cdot\uf_{>\up}(t^\uom\cdot \us)+\ep\cdot H(\uz,\us,\ep,t) =0.$
 In the case $\k=\C$ one takes the complex conjugation, $\overline{[\tilde\nabla \uf_\up|_\us]^T}.$

Finally, define the real (semi-algebraic) matrix $A\!:=\!\tilde\nabla \uf_\up|_\us\cdot\overline{[\tilde\nabla \uf_\up|_\us]^T}\!\in\! Mat_{c\times c}  .$
 Take its adjugate matrix, $A^\vee\!\cdot\! A\!=\!\one_{c\times c}\!\cdot\! det(A).$
 Multiply the last equation by $A^\vee$ to get:
 \beq\label{Eq.Foliation.of.semi-weighted.homogen.2}
 det\big[\tilde\nabla \uf_\up|_\us\cdot\overline{[\tilde\nabla \uf_\up|_\us]^T}\big] \cdot \uz
   +   A^\vee\cdot t^{-\de}\cdot t^{-\up}\cdot \uf_{>\up}(t^\uom\cdot \us)+\ep\cdot A^\vee\cdot H(\uz,\us,\ep,t) =0,
    \eeq
 for each  $\us\!\in\! Link[X_o]\sset \bS.$

 Recall the formula of Cauchy-Binet: $ det(B\cdot \overline{B^T})=\sum |det(\Box)|^2.$
  Here  $B\in Mat_{c\times N}$ is an arbitrary matrix, the sum goes over all the $c\times c$ blocks of $B$.
  In our case $det[\tilde\nabla \uf_\up|_\us\cdot\overline{[\tilde\nabla \uf_\up|_\us]^T]}=\sum |det(\Box)|^2,$ the sum goes over
    all the $c\times c$ blocks of $\tilde\nabla \uf_\up|_\us.$

Condition \eqref{Eq.Foliation.of.semi-weighted.homogen.2} should be satisfied for all the points $\us\in Link[X_o]\sset \bS .$ But all its ingredients are well defined at any point $\us\in \bS .$
 Thus we extend this condition to the whole standard sphere $\bS\sset \k^N $ as follows.
\\Take a polynomial $\tau\in \R[\us]$ satisfying: $\tau|_{Link[X_o]}=0$ and $\tau|_{\bS \smin Link[X_o]}>0.$
 E.g. one can take $\tau=\sum |f_{p_i}|^2|_{\bS }.$ Then equation \eqref{Eq.Foliation.of.semi-weighted.homogen.2} extends to:
 \beq\label{Eq.Inside.proof.of.foliations.Eq.to.resolve}
 \big(\tau+\sum |det(\Box)|^2\big) \cdot \uz
   +  A^\vee\cdot t^{-\de}\cdot t^{-\up}\cdot\uf_{>\up}(t^\uom\cdot \us)+\ep\cdot A^\vee\cdot H(\uz,\us,\ep,t) =0,
   \hspace{0.5cm} \forall\ \us\in \bS .
 \eeq
This is an implicit function equation on $\uz$, defined globally on $Ball_z\times Disc_\ep\times \bS \times [0,t_o].$
 The coefficient of the vector $\uz$
 is a non-negative polynomial, $\bS \to \R_{\ge0}.$ It vanishes exactly at those points of $Link[X_o]$ where  $rank[\tilde\nabla \uf_\up|_\us]<c.$
 Suppose a point $\us\in Link[X_o]$ satisfies: $\sum^{r_j}_{i=1}|s|^2_i=0,$ but $\sum^{r_{j+1}}_{i=1}|s|^2_i\neq0.$
  Then the condition $rank[\tilde\nabla \uf_\up|_\us]<r$ means: $rank[\di_{x_{r_j}+1}   \uf_\up|_\us,\dots,  \di_{x_N}   \uf_\up|_\us]<c.$
   Thus $s\in Sing[X_o\cap V(s_1\dots s_j)].$

Altogether, the coefficient $\tau+\sum |det(\Box)|^2$ vanishes exactly at the locus
 \beq
 Sing[X_o] \ \cup \ Sing[V(x_1,\dots,x_{r_1})\cap X_o] \ \cup  \ \cdots \ \cup \
 Sing[V(x_1,\dots, x_{r_k})\cap X_o].\vspace{-0.2cm}
 \eeq
And this is exactly the obstruction locus $\Si$ of equation \eqref{Eq.Def.of.Sigma}.

\

Altogether, \eqref{Eq.Inside.proof.of.foliations.Eq.to.resolve} is an implicit function equation. It is locally real-analytic in $\uz,t ,\ep,\us.$ %If the higher order terms, $\uf_{>\up},$ are algebraic power series, then the equation is an algebraic power series.

This equation is (uniquely) locally resolvable at each point $(\ep,\us)\in Disc\times (\bS \smin Link[\Si]).$
  For each fixed $(\ep,\us)$ its solution is an analytic  power series, $\uz(\ep,\us,t)\in \R\{t \},$
   resp. $\uz(\ep,\us,t)\in \R\bl t \br.$
   The radius of convergence
   can go to zero  as    $\us$ approaches $Link[\Si].$
 All these local solutions glue (by analytic uniqueness)  to a global solution. Therefore,  outside of an open (subanalytic) neighborhood
 $ Link[\Si]\sset \tcU \sset \bS,$ we get the global $\R$-analytic solution, $ \uz(\ep,\us,t)\in C^\om(Disc\times (\bS \smin \tcU )\times[0,t_o)) . $

If the higher order terms $\uf_{>\up}$ are Nash (i.e. $\R$-algebraic power series), then \eqref{Eq.Inside.proof.of.foliations.Eq.to.resolve} is an equation with algebraic power series. And thus the solutions are Nash,  $ \uz(\ep,\us,t)\in Nash(Disc\times (\bS \smin \tcU )\times[0,t_o)) . $

\item
For each point $(\ep,\us)\in Disc\times (\bS \smin Link[\Si]) $ we have constructed the perturbed arc
$\ga_{\ep,\us}(t)=(t^{\om_1}  s_1+t^{\om_1+\de} \ep\cdot h_1(\ep,\us,t) ,\dots,t^{\om_N} s_N+t^{\om_N+\de} \ep\cdot h_N(\ep,\us,t) ).$
 The  perturbations $\{h_\bullet\}$ satisfy:
$h_\bullet(\ep,\us,t)\in  C^\om(Disc\times (\bS\smin \tcU) \times[0,t_o)).$
%,\quad\quad ord_t h_\bullet(\ep,\us,t)\ge\de>0,\quad \quad
%\lim_{t\to 0}h_\bullet(\ep,\us,t)=0;
  If   $s_1=\cdots=s_{r_j}=0$   then  $h_1=\cdots=h_{r_j}=0.$

\

 \underline{The statements of part 2.}
\bee[\!\!\!\bf 2.i.\!]
\item  and {\bf 2.ii.} hold by the construction.
\setcounter{enumii}{2}
\item This follows immediately from 2.ii, by the ultrametric property of   $tord_{\k^N}(*,*).$

%The part $\Rrightarrow$ is immediate.
% For the part $\Lleftarrow$ we suppose that the arcs $\ga_\us,\ga_{\us'}$ are tangent. Then $\ga_\us,\ga_{\us'}\sset V(x_1,\dots,x_{r_j})\smin %V(x_1,\dots,x_{r_{j+1}}),$
%  for some (common) $0\le j\le k.$ First consider the case $j=0.$ One has:
%   $(s_1,\dots,s_{r_1})=(s'_1,\dots,s'_{r_1})\neq (0,\dots,0).$ Therefore  we have:
%   \beq
%   \ga_{\ep,\us}(t)= t (s_1,\dots,s_{r_1},0,\dots)+\ep\cdot t^{1+\frac{\de}{\om_1}}(\dots)+t^\frac{\om_2}{\om_1}(\dots) ,
%      \eeq
% \[\vspace{0.6cm}
%    \ga_{\ep',\us'}(t)=t (s_1,\dots,s_{r_1},0,\dots)+\ep'\cdot t^{1+\frac{\de}{\om_1}}(\dots)+t^\frac{\om_2}{\om_1}(\dots).
%\]
%   Hence the statement.
%The cases with $j>0$ (i.e. $s_1=\cdots=s_{r_j}=0$) go in the same way. Here we use the special form of $\uh$ in equation
 %\eqref{Eq.inside.proof.Def.of.h}.

   \item Assuming $s_1=\cdots=s_{r_j}=0, $ one has: $\ga_{\ep,\us}(t)=(0,\dots,0,s_{r_j+1}t^{\om_{r_j+1}}+\ep(\dots),\dots).$

\item If $\ga_\us\sset X_o\cap X_\ep\smin \Si$ for some $\ep\neq0$ and $\us\in \bS,$ then $\uf_\up(\ga_\us(t))=0=\uf_{>\up}(\ga_\us(t)).$ Then the unique solution of \eqref{Eq.Inside.proof.of.foliations.Eq.to.resolve} satisfies $z(\us)=0.$
\eee

  \

  \underline{The statements of part 3.}

     \bee[\!\!\!\!\!\!\!\bf 3.i.]
     \item If $\Si=\{o\}$ then $Link[\Si]=\empty,$ and equation \eqref{Eq.Inside.proof.of.foliations.Eq.to.resolve} is resolvable for each $(\ep,\us)\in Disc \times\bS $ on a corresponding segment $[0,t_o(\ep,\us)].$ By the compactness of $Disc \times\bS$ one has: $inf\{t_o(\ep,\us)\}>0.$ Thus the equation is resolvable globally,
           $\uz(\ep,\us,t)\in C^\om(Disc\times\bS\times[0,t_o])$ for some $0<t_o\ll1 .$

 If the higher order terms, $\uf_{>\up},$ are Nash power series, then we get the global Nash solution,
   $\uz(\ep,\us,t)\in Nash(Disc\times\bS\times[0,t_o)).$

\item To prove that $\{\ga_{\ep,\us}\}_\us$ is a foliation of $(\k^N,o)$, it is enough to rectify it or, rather, its strict transform. Define the map $\Psi_\ep:\bS\times[0,t_o]\to \bS\times[0,t_o]$ by
 $(s,t)\to (\frac{\us+\ep\cdot t^\de\cdot \uh}{\|\us+\ep\cdot t^\de\cdot \uh\|},t),$ where $\uh=\uh(s,t,\ep)$ is determined by our solution of \eqref{Eq.Inside.proof.of.foliations.Eq.to.resolve}.
  Note that $\uh$ is uniformly continuous on $\bS\times[0,t_o]\times Disc,$ thus bounded. Choose $0<t_o\ll1$ to get  $\|\ep\cdot t^\de\cdot \uh\|<\frac{1}{2}.$ Therefore $\Psi_\ep$ is a well-defined real-analytic map in $\us,t$.

We claim: $\Psi_\ep$ is ($\R$-analytically, globally) invertible for $0<t_o\ll1.$ It is enough to verify: for each fixed $0<t\le t_o$ the restriction  $\Psi_\ep|_t\circlearrowright\bS$ is an $\R$-analytic automorphism. Observe that $\Psi_\ep|_t$ is a deformation of the identity map
 $\Psi_\ep|_0=Id_\bS.$ Therefore for each fixed $\us_o\in \bS$ and $t\ll1$ the map-germ $\Psi_\ep|_t: (\bS,\us_o)\to (\bS,\Psi_\ep|_t(\us_o))$ is an isomorphism. And now invoke the compactness of $\bS,$ to get the needed uniform bound $0<t_o\ll1.$

Finally,    $\Psi_\ep(\tga_{\us})=\tga_{\ep,\us}.$ Therefore (for each $\ep$) the set of arcs $\{\tga_{\ep,\us}\}_\us$ is the image of the foliation
 $\{\tga_{\us}\}_\us$ under the   analytic isomorphism $\Psi_\ep.$ Hence $\{\tga_{\ep,\us}\}_\us$ is a (non-singular) foliation of
 $\bS\times[0,t_o],$ while $\{\ga_{\ep,\us}\}_\us$ is a (singular) foliation of
 $(\k^N,o).$

\eee

\

  \underline{The statements of part 4.}
 \bee[\!\!\!\!\!\!\!\bf 4.i.]
     \item Take a neighborhood $Link[\Si]\sset \tU\sset \bS.$
     Equation \eqref{Eq.Inside.proof.of.foliations.Eq.to.resolve} is resolvable locally for each $(\ep,\us)\in Disc\times(\bS\smin \tU).$ Thus, by compactness of $Disc\times(\bS\smin \tcU)$ we get the global $\R$-analytic solution
 $\uz(\ep,\us,t)\in C^\om(Disc\times(\bS\smin \tcU)\times[0,t_o)),$ for some $0<t_o\ll1.$

 And similarly in the Nash case.
\item Take the rectification of 3.ii,  $\Psi_\ep:(\bS\smin \tcU)\times[0,t_o]\to (\bS\smin Link[\Si])\!\times\![0,t_o].$
 It is defined for $0\!<\!t_o\!\ll\!1.$
  In particular, $\ga_{\ep,\us}\!\cap \!\ga_{\ep,\us'}\!\cap \!Ball_{r_o}\!=\!o$ for $s,s'\!\in \!\bS\!\smin \!\tcU$ with $s\!\neq \!s'.$

\item For each neighborhood $Link[\Si]\sset \tcU\sset \bS$ we get $0< t_o(\tcU)\ll1$ such that the properties 4.i, 4.ii hold on $(\bS\smin \tcU)\times[0,t_o(\tcU)).$ Thus there exists a (subanalytic) neighborhood $\tSi\sset \tcU(\tSi)\sset \bS\times[0,t_o)$ with the following properties, see \S\ref{Sec.Notations}:
 \beq
 Link_o[\tSi]\!=\!Link_o[\tcU(\tSi)], \quad
 Link[\tcU(\tSi)]\!=\!\tcU,\quad
 \text{the map $\Psi_\ep$ is defined on $\bS\!\times\![0,t_o]\!\smin\! \tcU(\tSi)$.}
 \eeq
Take the neighborhood   $\Si\sset \cU(\Si)\sset (\k^N,o)$ as the image of $\tcU(\tSi).$ It has the claimed properties.
\epr
\vspace{-0.4cm}
           \eee
\eee

\bex\label{Ex.Foliation.of.semi-weighted.homogen}
Let $X_o\!=\!V(f_p)\!\sset\! (\k^3,o)$ be a weighted-homogeneous surface germ with isolated singularity,
 $\om_1\!\le\! \om_2\!\le\!\om_3.$ Suppose $X_o$ does not contain the $\hx_3$-axis (of the top weight), while the intersection with the plane $V(x_1)\sset \k^3$ (of the bottom weight) is a reduced curve germ.
 Then equation \eqref{Eq.Def.of.Sigma} gives $\Si\!=\!\{o\}.$
 Therefore for any higher order perturbation $f_p\!+\!f_{>p},$ the weighted-homogeneous foliation $\{\ga_\us\}$ on $(\k^3,o)$ deforms into the foliation $\{\ga_{\ep,\us}\}$ compatible with $V(f_p+f_{>p}).$

When some of the weights coincide the assumptions are weakened:\vspace{-0.1cm}
\bee[$\bullet$]
\item Suppose $\om_1<\om_2=\om_3.$ Then just the intersection with the plane $V(x_1)$ should be a reduced curve germ.
\item Suppose $\om_1=\om_2<\om_3.$ Then it is enough to assume: $X_o\not\supset Span(\hx_3).$
\item If $\om_1=\om_2=\om_3$ then just the condition ``$X_o$ has an isolated singularity" suffices.
\eee
\eex

\beR\label{Rem.Deformation.of.Weight.Hom.Foliation.by.other.terms.}
 The perturbations in this lemma were by terms of higher order, $\uf_\up+\ep\cdot \uf_{>\up}.$
 One can allow also perturbations of the same order, $\uf_\up+\ep\cdot \uf_{\ge\up},$ thus $\de=0.$ But then the statement will be only for small $\ep,$ and without part $\Rrightarrow$ of $2.iii.$
 The proof is the same up to equation \eqref{Eq.Inside.proof.of.foliations.Eq.to.resolve}. And then one uses the assumption $|\ep|\ll1$ to resolve it.
\eeR

\bex\label{Ex.Briancon.Speder.}
  It is instructive to consider the Brian\c{c}on-Speder family, $f_\ep(x,y,z)=z^5+x^{15}+xy^7+\ep\cdot zy^6.$ It is weighted homogeneous with weights $(\frac{1}{15},\frac{2}{15},\frac{3}{15}).$ Equation \eqref{Eq.Def.of.Sigma} reads:
 \beq
 \Si=Sing[V(f_o)]\cup Sing[V(f_o,x)]\cup Sing[V(f_o,x,y)]\cup Sing[V(f_o,x,y,z)]\sset (\k^3,o).
 \eeq
 Here $Sing[V(f_o)]=\{o\}.$ The intersection $V(f_o,x)$ is non-reduced, thus $Sing[V(f_o,x)]=V(x,z).$
 Altogether,  $\Si=V(x,z)\sset (\C^3,o).$  Lemma \ref{Thm.semi-weighted.hom.foliation} gives: the weighted homogeneous foliation of $(\k^3,o)$ (which is compatible with $V(f_o)$) deforms to a foliation compatible with $V(f_o+\ep\cdot   zy^6)$ only outside of a hornic neighborhood of the $\hy$-axis.
\\
In fact this family of surface germ, $V(f_\ep)\!\!\sset\!\! (\C^3\!,o),$ is not inner-Lipschitz-trivial, \cite{Birbrair-Fernandes.Neumann.mu.constant.not.inner.Lip}. 
\mbox{Therefore the foliation cannot deform on the whole $(\C^3\!,o)$ compatibly with
    the surface $V(f_o\!+\!\ep\cdot zy^6).$}
\eex

\vspace{-0.4cm}
\section{Contact-trivialization of the family $\uf_\up+\ep\cdot \uf_{> \up}$}
\vspace{-0.1cm}
  The presentation of a function $g(\ux)$ in the weighted-polar coordinates(\S\ref{Sec.Notations}) is  $g(t^\uom\cdot \us).$
    Our functions/maps will be real-analytic off the origin, $C^\om(( \k^N\smin o,o)\times Disc),$ and globally analytic
       in  the weighted-polar coordinates,
  $g(t^{  \uom}\cdot \us)\in C^\om(\bS\times[0,t_o)\times Disc).$

Suppose  $V(\uf_\up)\sset (\k^N,o)$ is a weighted-homogeneous ICIS, and assume $\Si=o$ in \eqref{Eq.Def.of.Sigma}.

\bthe\label{Thm.Trivialization.of.Def.by.higher.order.terms}
The family $\uf_\up+\ep\cdot \uf_{> \up}$ is contact trivial, i.e. $\uf_\up+\ep\cdot \uf_{> \up}=U_\ep\cdot \uf_\up\circ \Psi_\ep,$ where:
 \bei
   \item $\Psi=\{\Psi_\ep\circlearrowright(\k^N,o)\}_\ep$ is a family of blow-analytic homeomorphisms, with $\Psi_o=Id_{(\k^N,o)};$
\item      $U=\{U_\ep \in GL(c)\}_\ep$ is a family of blow-analytic matrices, with $U_o=\one$.
\eei
   These $\Psi,U$ satisfy:
\bee[\!\!\bf 1.\!]
\item
  $\Psi\in C^\om(( \k^N\smin o,o)\times Disc) $ and $U \in GL\big(c,C^\om(( \k^N\smin o,o)\times Disc)\big).$
  \\
$\Psi_\ep$ is differentiable  at $o\in \k^n$ (for each $\ep $) with $\Psi'_\ep|_o=Id.$
  \\
  Their  weighted-polar presentations are globally real-analytic, i.e.
    \[
    \Psi (t^{ \uom}\cdot \us)\in C^\om(\bS\times[0,t_o)\times Disc) \quad \text{ and }\quad
 U (t^{ \uom}\cdot \us) \in GL(c,C^\om(\bS\times[0,t_o)\times Disc)).
\]
\item

Moreover, if the perturbations $ \uf_{>\up} $ \!are polynomials, then  $\Psi , U  $ are Nash on $(\k^N\!\smin \!o ,\!o)\!\times \!Disc,$    and their weighted-polar presentations are  Nash on $  \bS\times[0,t_o)\times Disc .$

\item (For each $\ep\in Disc$)
If $\de\!\ge \!\om_N-\om_2$ then $\Psi_\ep$ is bi-Lipschitz.
If
$\de\!>\! \om_N-\om_2$ then $\Psi_\ep$ is $C^1.$

\item (For each $\ep\in Disc$) If $\de\!\ge\! \om_N$ then $U_\ep$  is bi-Lipschitz.
 If $\de\!>\! \om_N$ then $U_\ep$  is   $C^1.$

\item In the hypersurface case ($c=1$) the family is $\cR$-trivial, i.e. one can choose $U_\ep=\one.$
\eee
\ethe
\bpr
In Step 1  we define $\Psi_\ep$ via the deformed foliation of Lemma \ref{Thm.semi-weighted.hom.foliation}. Using it we establish the equality of ideals in the weighted-polar coordinates, $(\uf_\up\circ\Psi_\ep)=(\uf_\up+\ep\cdot \uf_{> \up})\sset C^\om(\bS\times[0,t_o)).$ This implies the   contact triviality of the family. In Step 2 we verify the statements 1-5.
\bee[\!\!\!\bf Step 1.]
\item
 As $\Si\!=\!o,$ the foliation $\{\ga_\us\}$  of $(\k^N,o) $ deforms to $\ga_{\ep,\us},$   for each $\ep.$
  Thus we define the map $\Psi_\ep\!:(\k^N,o)\!\to\! (\k^N,o)$ by $\ga_\us(t)\!\to\! \ga_{\ep,\us}(t).$
  Explicitly, $\Psi_\ep: x_\us(t)\to x_\us(t)\!+\!\ep\!\cdot\! t^\de\!\cdot\! t^{\uom}\!\cdot \! \uh(\us,t,\ep),$ where $\uh$ is given by the Ansatz  \eqref{Eq.inside.proof.Def.of.h}, while $\uz$ is determined by equation \eqref{Eq.Inside.proof.of.foliations.Eq.to.resolve}.
 Thus $\Psi_\ep$ is a family of (subanalytic) homeomorphisms, and their weighted-polar presentation is  analytic/Nash.
\\
The map of weighted-polar coordinates, $\bS\times\R_{>0}\to \k^N\smin o ,$ is a real-algebraic isomorphism.
 (Note that $gcd(\om_1,\dots,\om_N)\!=\!1.$)
  Thus $\Psi_\ep$ is obviously real-analytic/Nash over
 $(\k^N\smin o,o).$

  By the construction $\Psi$ trivializes the family of zero-sets,
\beq
V(\uf_\up+\ep\cdot \uf_{> \up})=\Psi_\ep (V(\uf_\up))=V(\uf_\up\circ \Psi^{-1}_\ep) \sset (\k^N,o).
\vspace{-0.2cm}\eeq
  We promote this to the contact trivialization of functions.
 \bee[\!\!\!\!$\bullet$]
 \item
 (The hypersurface case, $c\!=\!1$.) Take the weighted-polar presentation $f_p(t^\uom\!\cdot  \!\us )\!+\!\ep\!\cdot \!f_{> p}(t^\uom\!\cdot  \!\us )\!\in \!C^\om(\bS\!\times\![0,t_o)).$ We have the $\R$-analytic trivialization of the family of zero sets, $V(f_p(t^\uom\cdot  \us)+\ep\cdot f_{> p}(t^\uom\cdot  \us))=V(f_p(\Psi^{-1}_\ep(t^\uom\cdot  \us))).$ Therefore locally at each point of $\bS\times[0,t_o)$ the family of functions
  $f_p(t^\uom\cdot  \us )+\ep\cdot f_{> p}(t^\uom\cdot  \us )$ is $\R$-analytically contact trivial.
 (See e.g. pg.113, Proposition 4.6 in \cite{Mond-Nuno}.)

 Therefore the ratio $\tu(t,\us):=\frac{f_p(t^\uom\cdot  \us )+\ep\cdot f_{> p}(t^\uom\cdot  \us )}{f_p(\Psi^{-1}_\ep(t^\uom\cdot  \us))}$
  is locally analytic and invertible at each point of $\bS\times[0,t_o).$ Thus (by analytic uniqueness) we   get the globally analytic functions
   $\frac{1}{\tu},\tu\in C^\om(\bS\times[0,t_o)).$

   Finally we define the ``scaling factor" $u:(\k^N,o)\to\k^1$ via $u(t^\uom\cdot  \us)=\tu(t,\us)$ to get the contact trivialization
    $f_p+\ep\cdot f_{> p}=u\cdot f_p\circ\Psi^{-1}_\ep.$

 \item (The general case, $c\ge1$.) As before we have the  $\R$-analytic trivialization of the family of zero sets, $V(\uf_\up(t^\uom\cdot  \us)+\ep\cdot \uf_{> \up}(t^\uom\cdot  \us))=V(\uf_\up(\Psi^{-1}_\ep(t^\uom\cdot  \us)))\sset \bS\times[0,t_o).$
 As before (by Mather's theorem), at each point  of $\bS\times[0,t_o)$ the family of functions
  $\uf_\up(t^\uom\cdot  \us )+\ep\cdot \uf_{> \up}(t^\uom\cdot  \us )$ is locally $\R$-analytically contact trivial.
  Thus the germs of the ideals $(\uf_\up(t^\uom\cdot  \us )+\ep\cdot \uf_{> \up}(t^\uom\cdot  \us )),(\uf_\up(\Psi^{-1}_\ep(t^\uom\cdot  \us)))\sset C^\om(\bS\times[0,t_o))$ coincide locally at each point of $\bS\times[0,t_o).$
   We want to deduce: these ideals coincide globally.

 Consider the quotient module
 \beq
 M:=\quot{( \uf_\up(\us)+t^{-\up}\cdot \ep\cdot \uf_{>\up}(t^{\uom}\cdot\us))+(\uf_\up(\us))}{(\uf_\up(\us))}\in mod\text{-}C^\om(\bS\times(-t_o,t_o)).\vspace{-0.2cm}
 \eeq
  We want to prove  $M=0.$ This module is finitely generated, and its fibres vanish at all the points of
 $\bS\times(-t_o,t_o).$
  \bee[\!\!\!$-$]
  \item Working over $\C$ we would immediately conclude: $M=0.$ (See e.g. pg. 46 of \cite{Atiyah-Macdonald}.) Indeed, this is equivalent to   
  $Supp M=\empty.$ Thus it is enough to verify  that the localizations of $M$ at all the maximal ideals vanish. 
   Thus ($M$ being a f.g. module) it is enough to verify that the fibres of $M$ at all the f.g. maximal ideals vanish. 
  And finally,  the f.g. maximal ideals correspond to the closed points.
 \item   Over $\R$ this is not immediate, as many maximal ideals do not correspond to the real closed points of $\bS\times(-t_o,t_o).$
  Therefore we observe: $Supp M$ is a (Zariski) closed subset of a scheme whose real points are the points of $\bS\times(-t_o,t_o).$
  And $Supp(M)$ contains no point of $\bS\times(-t_o,t_o).$ Thus $Supp(M)$ is disjoint to a neighborhood of the real locus $\bS\times(-t_o,t_o).$
   And thus $M$ vanishes along this locus.
\eee
   Altogether we get the global equality of ideals:  $(\uf_\up(\us))\!=\!( \uf_\up(\us)\!+\!t^{-\up}\cdot \ep\cdot \uf_{>\up}(t^{\uom}\cdot\us)).$ Which means: $\uf_\up(\us)\!+\!t^{-\up}\cdot \ep\cdot \uf_{>\up}(t^{\uom}\cdot\us)= U\cdot \uf_\up(\us),$
    for a matrix $ U(t^{\uom}\cdot \us,\ep)\in GL(c,C^\om\big(\bS\times(-t_o,t_o)\times Disc)\big).$
 \eee

 \item
{\bf Statement 1.} We should only verify the differentiability of $\Psi_\ep\circlearrowright(\k^N,o)$ at the origin.
  It is enough to present
  $\uh(\ux)=B(s,t)\cdot \ux,$ where the matrix $B\in Mat_{N\times N}$ is bounded.
 Write the entries of  $t^\uom\cdot \uh:$\vspace{-0.4cm}
 \beq
t^{\om_j} h_j(\ep,\us,t)=(\sum^{r_l}_{i=1} |s^2_i|)\cdot t^{\om_j}\cdot \Big[\di_{x_j}f_{p_1}|_\us,\dots,\di_{x_j}f_{p_c}|_\us\Big]\cdot A^\vee\cdot \bbm t^{-p_1-\de}\cdot f_{>p_1}(t^{\uom}\cdot\us)+(h.o.t.)\\\dots\\t^{-p_c-\de}\cdot f_{>p_c}(t^{\uom}\cdot\us)+(h.o.t.)\ebm.
 \eeq
Here   $r_l\ge j$ and $\om_j=\om_{r_l},$ and $s_i t^{\om_j}=x_i\cdot t^{\om_j-\om_i}, $ with $\om_j\ge \om_i.$
  Therefore we can present: $t^{\om_j} h_j(\us,t,\ep)=\sum^{r_1}_{i=1} x_i\cdot b_{i,j}(\us,t,\ep).$  And one has $ord_\om(h.o.t.)> ord_\om(t^{-p_k-\de}\cdot f_{>p_k})\ge 0.$ Thus  the coefficients $b_{i,j}(\us,t,\ep)$ are bounded.

 With this presentation we have: $\Psi_\ep(\ux)-\ux=\ep\cdot t^\de\cdot B(s,t)\cdot \ux=o(\|\ux\|).$ Hence $\Psi_\ep$  is differentiable, and
 $\Psi'|_o=Id.$

{\bf Statement 2} is immediate.

{\bf Statement 3.} We  bound  the  derivatives $\frac{\di (t^{\de+\om_i}h_i)}{\di x_j}.$
 It is enough to verify that the polar presentation of these derivatives is bounded (resp. tend to $0$) when $t\to 0.$
  By compactness of $\bS$ it is enough to verify this locally at each point.
  Finally, it is enough to consider the case  $\k=\R.$ (For $\k=\C$ one passes to the real-polar coordinates, $s_1,\dots,s_{2N}$.)

It is enough to work in the chart $s_N\neq0$, where the coordinates are $t,s_1,\dots,s_{N-1}.$
  (Then $s_N$ is determined by $\sum^N |s_i|^2=1$.) Indeed, this chart is dense on $\bS.$

We use the chain rule,
$\frac{\di(t^{\de+\om_1}h_1,\dots,t^{\de+\om_N}h_N)}{\di(x_1\dots x_N)}=\frac{\di(t^{\de+\om_1}h_1,\dots,t^{\de+\om_N}h_N)}{\di(t,s_1\dots s_{N-1})}
\cdot \frac{\di(t,s_1\dots s_{N-1})}{\di(x_1\dots x_N)}.$
 To compute the matrix $\frac{\di(t,s_1\dots s_{N-1})}{\di(x_1\dots x_N)}$ we write its inverse:\vspace{-0.2cm}
\beq
\frac{\di(x_1\dots x_N)}{\di(t,s_1\dots s_{N-1})}=
\bbm
 s_1 t^{\om_1-1}&t^{\om_1}&0&\dots \\
 s_2 t^{\om_2-1}&0&t^{\om_2}&0
\\\dots&\dots&0&t^{\om_{N-1}}
\\ s_N t^{\om_N-1}&\frac{\di s_N}{\di s_1}t^{\om_N}&\dots& \frac{\di s_N}{\di s_{N-1}}t^{\om_N}
 \ebm.
\eeq
 Therefore $det[\frac{\di(x_1\dots x_N)}{\di(t,s_1\dots s_{N-1})}]\stackrel{*}{=}t^{\sum \om_i-1}(s_1 \frac{\di s_N}{\di s_1}+\cdots+ s_{N-1} \frac{\di s_{N}}{\di s_{N-1}}+s_N)=s_N\cdot t^{\sum \om_i-1}.$
  Here $*$ is obtained by extracting  $t^{\om_i}$ from $i$'th row and $\frac{1}{t}$ from the 1'st column.

  The adjugate matrix is $[\frac{\di(x_1\dots x_N)}{\di(t,s_1\dots s_{N-1})}]^\vee\in t^{\sum \om_i-1-\om_N}\cdot Mat_{N\times N}(\R\{\us,t\}).$ The first tow of this matrix belongs to the ideal $(t^{\sum \om_i-\om_N}).$
  Thus $\frac{\di(t,s_1\dots s_{N-1})}{\di(x_1\dots x_N)}\in t^{-\om_N}\cdot Mat_{N\times N}(\R\{\us,t\}).$ And the first row of this matrix belongs to the ideal $t^{1-\om_N}.$

 Altogether we get:
$ \frac{\di(t^{\de+\om_1}h_1,\dots,t^{\de+\om_N}h_N)}{\di(x_1\dots x_N)}\in t^{\de+\om_1-\om_N}\cdot Mat_{N\times N}(\R\{\us,t\}).$
  Thus, for $\de\ge \om_N-\om_1,$ the derivative of $\Psi_\ep$  is bounded on $(\k^N,o).$ Hence $\Psi_\ep$ is Lipschitz.
   If  $\de> \om_N-\om_1,$ the derivative of $\Psi_\ep$  is continuous on $(\k^N,o).$ Hence $\Psi_\ep$ is $C^1.$
  The corresponding statements for $\Psi^{-1}_\ep$ are obtained in the same way.

In the case $\om_1=\om_2$ we have proved statement 3. If $\om_1<\om_2$  then $h_1$ is divisible by $s_1.$ Therefore, once we have found the solutions
 $h_1(t,\us,\ep),\dots, h_N(t,\us,\ep)$ we can rescale $t,$ by $s_1\cdot \tt^{\om_1}=s_1\cdot t^{\om_1}(1+\frac{h_1}{s_1}),$ to get: $h_1=0.$
  Then we get the improved bound:
  $ \frac{\di(0,t^{\de+\om_2}h_2,\dots,t^{\de+\om_N}h_N)}{\di(x_1\dots x_N)}\in t^{\de+\om_2-\om_N}\cdot Mat_{N\times N}(\R\{\us,t\}).$ This gives the statement 3.
\\{\bf Statement 4} follows by observing: $U-\one\in (t^\de)$ and  $\frac{\di(t,s_1\dots s_{N-1})}{\di(x_1\dots x_N)}\in t^{-\om_N}\cdot Mat_{N\times N}(\R\{\us,t\}).$
\\
{\bf Statement 5.} For $c=1$   we use the weighted homogeneity of $f_p$ to push $U_\ep$ inside $f_p.$ The corresponding adjusted homeomorphism $\Psi_\ep$
  still has all the claimed properties.
\epr\vspace{-0.5cm}
\eee

\vspace{-0.2cm}
\bex\label{Ex.Trivialiations}
(Comparison  to the known statements on contact/right Lipschitz trivialization.) Take a map $\uf_\up:\R^N\to \R^c$ such that $\uf^{-1}_\up(o)$ is an ICIS.\vspace{-0.2cm}
\bee[\!\!\bf i.\!]
\item
\bei
\item   \cite[Theorem 3.3]{Fernandes-Ruas},    for $f_p:(\R^N,o)\to(\R,o)$:
\\If $\de\ge \om_N-\om_1$ then the family $f_p+\ep\cdot f_{>p}$ is Lipschitz-right trivial.
 \item  \cite[Proposition 2.2]{Ruas-Saia.97},   for $f_p:(\R^N,o)\to(\R,o)$:
 \\If $\de\ge \om_N-\om_1+1$ then the family $f_p+\ep\cdot f_{>p}$ is $C^1$-right trivial.
\item \cite[Theorem 2.3]{Costa-Saia-Soares},  for $\uf_\up:(\R^N,o)\to(\R^c,o)$:
 \\If $\de\ge \om_N-\om_1  $   then the family $\uf_\up+\ep\cdot \uf_{>\up}$ is Lipschitz-contact trivial.
  \eei
Theorem \ref{Thm.Trivialization.of.Def.by.higher.order.terms}   slightly improves these bounds.
 What is more important, the  trivializations in these papers are obtained by vector field integration. It is not clear whether these trivializations
 are subanalytic or differentiable, not speaking about the blow-analyticity.

On the other hand, \cite[Theorem 0.2]{Fukui.Paunescu.00} prove (for $f_p:(\R^N,o)\to(\R,o)$): the family $f_p+\ep\cdot f_{>p}$ is right trivializable by a homeomorphism whose strict transform  under the weighted-homogeneous blowup is real-analytic. But (it seems)  they do not have any control on the Lipschitz/$C^1$-properties of  that trivialization.

  Finally,  the trivializations of all these mentioned papers, were obtained by vector field integration.
  Thus the   Nash'ness on $(\k^N\smin o,o)$ and on $\bC\times[0,t_o)$ were beyond any hope.

\item Let $\om_1=\cdots=\om_N ,$ and assume $\Si=o.$  We get: the family $\uf_\up+\ep\cdot \uf_{>\up}$ is contact-$C^1$-trivial.
 In the case of one function ($c=1$) the family is   right-$C^1$-trivial.
  This is well known, see e.g. \cite{Takens} and \cite{Costa-Saia-Soares}.
\eee
\medskip
\eex
In theorem \ref{Thm.Trivialization.of.Def.by.higher.order.terms} the perturbations $\uf_{>\up}$ were by terms of higher $t$-order.
 We can also  perturbations of the same order, $\uf_{\up,\ep}+\ep\cdot \uf_{\ge\up},$ then getting the weaker version.
\bcor\label{Thm.Trivialization.for.defs.by.same.order.terms}
Suppose $\Si(\uf_{\up,\ep})=o$ for each $\ep\in Disc.$ Then the family $\uf_{\up,\ep}+\ep\cdot \uf_{\ge\up}$ is contact-trivial,
 $\uf_{\up,\ep}+\ep\cdot \uf_{\ge\up}=U_\ep \cdot \uf_\up\circ\Psi_\ep,$ where $U_\ep,\Psi_\ep$ satisfy properties 1,2,5 of theorem \ref{Thm.Trivialization.of.Def.by.higher.order.terms}, except for the differentiability of $\Psi_\ep$, but not necessarily properties 3,4.
\vspace{-0.2cm}\ecor
\bpr
 For the case $ |\ep|\ll1$  the proof is the same as in theorem \ref{Thm.Trivialization.of.Def.by.higher.order.terms}. The homeomorphism $\Psi_\ep$ is again defined via the deformed foliation, see Remark \ref{Rem.Deformation.of.Weight.Hom.Foliation.by.other.terms.}.

In the general (non-local) case, $\ep\in Disc,$ one takes a finite cover of the compact disc by small balls. For each such small ball we get
 the unique coordinate change $\Psi_\ep$ whose weighted-polar presentation is $\R$-analytic. These coordinate changes glue by analytic uniqueness.
  (Note that $Disc$ is simply-connected, thus we have no monodromy obstructions.) Thus we get the global family of homeomorphisms $\Psi_\ep\in C^\om(Disc\times \bS\times[0,t_o)).$ By the same arguments as in theorem this leads to the global equality of ideals
   $(\uf_\up\circ\Psi_\ep)=((\uf_\up+\ep\cdot \uf_{>\up})\circ\Psi_\ep)\sset C^\om(\bS\times[0,t_o)).$ This gives the statement.
\epr

\end{document}